\documentclass{amsart}

\usepackage{graphicx}
\usepackage{enumerate}
\usepackage{amssymb}

\newtheorem{theorem}{Theorem}[section]
\newtheorem{lemma}[theorem]{Lemma}

\newtheorem{mtheorem}{Theorem}

\newtheorem{mcorollary}[mtheorem]{Corollary}

\theoremstyle{definition}
\newtheorem{definition}[theorem]{Definition}
\newtheorem{remark}[theorem]{Remark}

\numberwithin{equation}{section}
\numberwithin{figure}{section}

\begin{document}

%%%%%%%%%%%%%
\title[Generic cubic homoclinic tangencies for H\'enon maps]{Existence of generic cubic homoclinic tangencies for H\'enon maps}

\author{Shin Kiriki}
\address{Department of Mathematics, Kyoto University of Education, 
1 Fukakusa-Fujinomori, Fushimi-ku, Kyoto, 612-8522, JAPAN}
\email{skiriki@kyokyo-u.ac.jp}

\author{Teruhiko Soma}
\address{Department of Mathematics and Information Sciences,
Tokyo Metropolitan University,
Minami-Ohsawa 1-1, Hachioji, Tokyo 192-0397, JAPAN}
\email{tsoma@tmu.ac.jp}

\date{\today}

\maketitle
%%%%%%%%%%%%%

\begin{abstract}
In this paper, we show that the H\'enon map $\varphi_{a,b}$ has a generically unfolding cubic tangency for some $(a,b)$ arbitrarily close to $(-2,0)$ by applying results of Gonchenko-Shilnikov-Turaev \cite{GST1}\,--\,\cite{GST7}.
Combining this fact with theorems in Kiriki-Soma \cite{KS}, 
one can observe the new phenomena in the H\'enon family, appearance of persistent antimonotonic tangencies and cubic polynomial-like strange attractors.
\end{abstract}

\setcounter{section}{-1}

\section{Introduction}

Let $f$ be a 2-dimensional diffeomorphism with a cubic homoclinic tangency associated with a dissipative saddle fixed point.
In \cite{KS}, we showed that, if the cubic tangency unfolds generically in a two-parameter family 
$\{f_{\mu,\nu}\}$ with $f=f_{0,0}$, then  the family exhibits cubic dynamics, the existence of persistent 
antimonotonic tangencies and cubic polynomial-like strange attractors.
The dynamics is quite different from quadratic dynamics inherited from 
one-dimensional maps.
If a diffeomorphism with a cubic homoclinic tangency as above is given, then it is not hard to find such 
a two-parameter family.
However, the converse is not obvious.

In this paper, we consider the case that the two-parameter family is  
the H\'enon family introduced by \cite{He}.
\emph{H\'enon maps}
\begin{equation}\label{original}
f_{a,b}(x,y) = (1+y -ax^2,bx)\quad (b\neq 0)
\end{equation}
are one of the most important models of two-dimensional diffeomorphisms in modern chaotic dynamical systems.
When $b$ is sufficiently small, it seems that the dynamics of $f_{a,b}$ is similar to that of the 
quadratic map $f_a(x)=1-ax^2$ on $\mathbb{R}$.
In fact, for any sufficiently small $b>0$, Benedicks and Carleson \cite{BC2} found a positive Lebesgue measure subset $J_b$ of $a$-values near $2$ such that $f_{a,b}$ has a quadratic polynomial-like strange attractor if $a\in J_b$ by generalizing the result on the quadratic maps \cite{BC}, see \cite{LV, MV} for related results.
We refer to \cite{BV,BY,WY} for ergodic properties on such strange attractors, see also \cite{BDV} for comprehensive references on concerned topics.

For convenience in our arguments, we adopt the following formula of the H\'enon map: 
$$
\varphi_{a,b}(x,y)=(y, a-bx+y^2) 
$$
which is obtained from the classical formula (\ref{original}) by the 
 reparametrization $(a,b)\mapsto(-a,-b)$ and the coordinate change $(x,y)\mapsto (-ab^{-1}y, -ax)$.
The fixed points of $\varphi_{a,b}$ are $p_{a,b}^\pm=(y_{a,b}^\pm,y_{a,b}^\pm)\in \mathbb{R}^2$ with $y_{a,b}^\pm=\bigl(1+b\pm\sqrt{(1+b)^2-4a}\bigr)/2$.
Note that $p^+_{a,b}$ and $p^-_{a,b}$ converge respectively to the points $(2,2)$ and $(-1,-1)$ 
in $\mathbb{R}^2$ as $(a,b)\rightarrow(-2,0)$. 

The following is our main theorem.

\begin{mtheorem} \label{main_a}
There exists $(a^*,b^*)$ with $b^*>0$ and arbitrarily close to $(-2,0)$  
such that the H\'enon map $\varphi_{a^*,b^*}$ has a cubic homoclinic tangency associated with $p_{a^*,b^*}^{+}$ which unfolds generically with respect to $\{\varphi_{a,b}\}$.
\end{mtheorem}

Moreover, the tangency can be supposed to be of type I in the sense of Definition \ref{dfn_cub_type}.
This fact is used to prove Corollary \ref{cor_b}.

\begin{figure}[htb]
\begin{center}
\scalebox{0.8}{\includegraphics[clip]{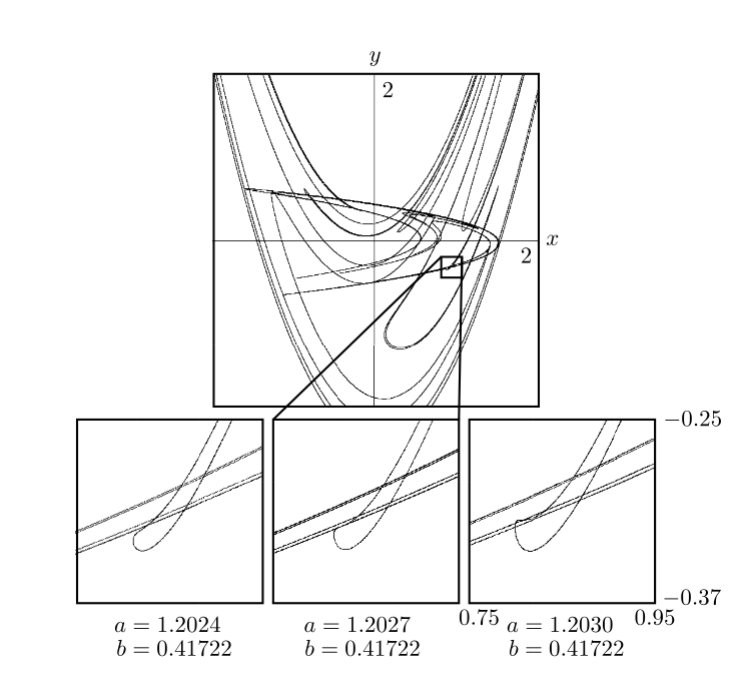}}
\caption{}
\label{fig_01}
\end{center}
\end{figure}

Carvalho \cite[p.\ 769]{Ca} presents as numerical results a supporting evidence of the existence of generically unfolding cubic homoclinic tangencies in the H\'enon family at parameters $(a,b)$ near $(1.203, 0.417)$ and $(1.095,0.388)$.
Figure \ref{fig_01} illustrates the stable and unstable manifolds of H\'enon maps $f_{a,b}$ with $(a,b)$ around $(1.2027, 0.41722)$, which are depicted by using the software \textit{Janet}\footnote{Available from C. Knudsen's home page: \tt{http://dcwww.fys.dtu.dk/$\sim$carsten/}} produced by Knudsen et al.
%%%%%%%%%%%%%%%%%%%%%%%%%
%%%%%%%%%%%%%%%%%%%%%%
However, as far as the authors know, any strict proof of the existence of such tangencies has not been obtained.
Our proof is based on the mechanism introduced by Gonchenko, Shilnikov and Turaev \cite{GST1}--\cite{GST7} 
which produces a cubic homoclinic tangency from a pair of quadratic heteroclinic tangencies cyclically associated with a pair of saddle points.

As for Anosov diffeomorphisms, 
Bonatti, D\'iaz and Vuillemin \cite{BDVu} detect a codimension-three submanifold 
$\mathcal{C}$ of $\mathrm{Diff}^3(T^2)$ contained in the boundary of Anosov 
diffeomorphisms such that each element of $\mathcal{C}$ has a cubic heteroclinic tangency.
See also \cite{Ca} for related results.

Here we propose the following question  
naturally arising from Theorem \ref{main_a} and Gonchenko-Shilnikov-Turaev [14, Theorem 1]. 

\newtheorem{question}{Open Question}
\renewcommand{\thequestion}{\!}
\begin{question}
Does the original H\'enon family (\ref{original})  have a homoclinic tangency of every order?
\end{question}

The following corollary is obtained immediately from our Theorem \ref{main_a} together with Theorem \ref{T_KS} in Subsection \ref{SS_renormal}, which is a stronger version of Kiriki-Soma \cite[Theorems A 
and B]{KS}. 
From this corollary, we know that the H\'enon family has cubic dynamics as well as quadratic dynamics.

\begin{mcorollary}\label{cor_b}
There exist subsets $\mathcal{O}$ and $\mathcal{Z}$ in the $ab$-space such that the H\'enon subfamilies $\{\varphi_{a,b}\}_{(a,b)\in \mathcal{O}}$ and  $\{\varphi_{a,b}\}_{(a,b)\in \mathcal{Z}}$ satisfy the following conditions.
\begin{enumerate}[\rm (i)]
\item 
$\mathcal{O}$ is an open set with $\mathrm{Cl}(\mathcal{O})\ni (-2,0)$.
For any $(a,b)\in \mathcal{O}$ and a sufficiently small $\varepsilon>0$, there exists a regular curve $c:(-\varepsilon,\varepsilon)\rightarrow \mathcal{O}$ with $c(0)=(a,b)$ such that the one-parameter family $\{\varphi_{c(t)}\}$ exhibits persistent antimonotonic tangencies.
\item 
For any open neighborhood $U$ of $(-2,0)$ in the $ab$-space, $\mathcal{Z}\cap U$ has positive $2$-dimensional Lebesgue measure.
For any $(a,b)\in \mathcal{Z}$, there exists an integer $n+N>0$ such that $\varphi_{a,b}^{n+N}$ exhibits a cubic polynomial-like strange attractor of type I supported by an SRB measure.
\end{enumerate}
\end{mcorollary}

\begin{figure}[htb]
\begin{center}
\scalebox{0.8}{\includegraphics[clip]{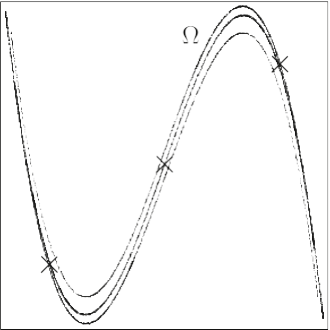}}
\caption{}
\label{fig_02}
\end{center}
\end{figure}

An invariant set $\Omega$ of a $2$-dimensional diffeomorphism $\psi$ is called a \emph{strange attractor} if (a) there exists a saddle point $p\in \Omega$ such that the unstable manifold $W^u(p)$ has dimension $1$ and $\mathrm{Cl}(W^u(p))=\Omega$, (b) there exists an open neighborhood $U$ of $\Omega$ such that $\{\psi^n(U)\}_{n=1}^\infty$ is a decreasing sequence with $\Omega=\bigcap_{n=1}^\infty \psi^n(U)$, and (c) there exists a point $z_0\in \Omega$ whose positive orbit is dense in $\Omega$ and a non-zero vector $\boldsymbol{v}_0\in T_{z_0}(\mathbb{R}^2)$ with $\Vert D\psi^n(z_0)(\boldsymbol{v}_0)\Vert \geq e^{cn}\Vert \boldsymbol{v}_0\Vert$ for any integer $n\geq 0$ and some constant $c>0$.
Besides $\Omega$ is called a \textit{cubic polynomial-like} strange attractor \emph{of type I} if there exists an integer $m>0$ such that the dynamics of $\psi^m$ on $\Omega$ is close 
(up to scale) to that of the one-dimensional map $x\mapsto -x^3+ax$ with $a\in (3\sqrt{3}/2,3)$ and has three saddle fixed points, see Figure \ref{fig_02}.
%%%%%%%%%%%%%%%%%%%%%%%%%
%%%%%%%%%%%%%%%%%%%%%%
An \textit{SRB measure} on $\Omega$ means  a $\psi$-invariant Borel probability measure  which is ergodic, has a compact support and 
has absolutely continuous conditional measures on unstable manifolds.

The theorems in \cite{KS} were proved by invoking a renormalization of two-dimensional $C^\infty$ diffeomorphisms near a cubic homoclinic tangency (\cite[Lemma 2.1]{KS}) together with some results of Wang-Young \cite{WY}.
In this paper, in stead of \cite[Lemma 2.1]{KS}, we apply  Lemma \ref{G_Lemma_6} in Subsection \ref{SS_renormal}, 
which is a special case of the renormalization presented by 
Gonchenko et al \cite[Lemma 6]{GST7}.
By this replacement, we have an improved theorem, Theorem \ref{T_KS}, which does not need any longer the linearizability condition on $\varphi_{a,b}$ near $p_{a,b}^-$ used in the 
proofs of the theorems in \cite{KS} and moreover it works even in the $C^r$ category 
with $4\leq r\leq \omega$.
Throughout this paper, we regard H\'enon maps as such $C^r$ diffeomorphisms.
So Corollary \ref{cor_b} will hold for any two-parameter 
family in $\mathrm{Diff}^r(\mathbb{R}^2)$ sufficiently close to the H\'enon family.

Kan, Ko\c{c}ak and Yorke \cite{KKY} showed that, for any one-parameter family 
$\{f_\mu\}$ of two-dimensional diffeomorphisms such that $f=f_0$ has a generically unfolding 
quadratic homoclinic tangency associated with a dissipative saddle point, 
infinitely many $f_\mu$ have both contact-making 
and contact-breaking tangencies simultaneously, called \emph{antimonotonic tangencies}.
This is contrast to the fact by Milnor-Thurston \cite{MT} that the one-dimensional quadratic maps $x\mapsto 1-ax^2$ can have only 
orbit-creation values and no orbit-annihilation values, that is, periodic orbits are created monotonically 
as the parameter $a$ increases.
Our persistent antimonotonic tangencies in Corollary \ref{cor_b} are \emph{cubically related} in the sense of \cite[Subsection 4.2]{KS}, which have properties different from the antimonotonic tangencies given in \cite{KKY}.

This paper is organized as follows.
Section \ref{S_Pre} presents fundamental notations and definitions needed in later sections.
In Section \ref{S_GST}, results of Gonchenko-Shilnikov-Turaev \cite{GST1}--\cite{GST7} are revisited.
In particular, we will review their results on sufficient conditions for the existence of a cubic homoclinic tangency in a two-parameter family of two-dimensional diffeomorphisms. 
In Section \ref{S_3}, we will prove Theorem \ref{main_a}.
An outline of our proof is given in Subsection \ref{S_outline_b}.

\section{Generically unfolding tangencies}\label{S_Pre}

In this section, we will review some properties of generically unfolding quadratic and cubic tangencies 
for 2-dimensional diffeomorphisms.
Though they are more or less known results, we present them in forms suitable to our arguments.

\subsection{Generically unfolding quadratic tangencies}\label{S_gqt}

A diffeomorphism $\psi$ on $\mathbb{R}^2$ has a \textit{transverse point} $r$ associated with saddle fixed points $p_1, p_2$ if
\begin{itemize}
\item $r\in W^{u}(p_{1})\cap W^{s}(p_{2})\setminus\{p_{1},p_{2}\}$,
\item $\dim(T_{r}W^{u}(p_{1})+T_{r}W^{s}(p_{2}))=2$.
\end{itemize}
We also say that $\psi$ has a \textit{tangency} $q$ of order $n$ 
associated with saddle fixed points $p_{1}, p_{2}$ if it satisfies the following conditions.
\begin{itemize}
\item $q\in W^{u}(p_{1})\cap W^{s}(p_{2})\setminus\{p_{1},p_{2} \}$.
\item $\dim(T_{q}W^{u}(p_{1})+T_{q}W^{s}(p_{2}))=1$.
\item There exists a local $C^{n+1}$ coordinate $(x, y)$ in a neighborhood of $q$ such that 
$q = (0,0)$, 
$\{(x, y) ; y=0\}\subset W^{s}(p_{1})$ and $\{(x, y); y=\alpha(x)\}\subset W^{u}(p_{2})$, 
where $\alpha$ is a $C^{n+1}$-function satisfying
\begin{equation}\label{eqn_tangency}
\alpha(0)=\alpha^{\prime}(0)=\cdots=\alpha^{(n)}(0)=0\quad\mbox{and}\quad
\alpha^{(n+1)}(0)\neq 0.
\end{equation}
\end{itemize}
In the case when $p_{1}=p_{2}$, the transverse point or the tangency is called to be \textit{homoclinic}, and otherwise \textit{heteroclinic}. 
The definition of a tangency of order $n$ is independent of the choice of a local $C^{n+1}$ coordinate satisfying the condition as above.
Usually, the first order tangency is called \textit{quadratic}, and the second order is \textit{cubic}.

One parameter family $\{Q_\mu\}$ with $Q_\mu\subset \mathbb{R}^2$ is called a set of $C^r$ 
\emph{continuations} (for short \emph{continuations}) of curves in $\mathbb{R}^2$ if there exist $C^r$ embeddings $f_\mu:I\to 
\mathbb{R}^2$ of a fixed interval $I$ which $C^r$ depend on $\mu$ and satisfy 
$f_\mu(I)=Q_\mu$.
A set of continuations $\{q_\mu\}$ of points in $\mathbb{R}^2$ is defined similarly.

Let $\{\psi_\mu\}_{\mu\in J}$ be a one-parameter family in $\mathrm{Diff}^r(\mathbb{R}^2)$ $(r\geq 2)$ such that the parameter space $J$ is an interval, and $p_{1,\mu},p_{2,\mu}$ (possibly $p_{1,\mu}=p_{2,\mu}$) continuations of saddle fixed points of $\psi_\mu$ such that $W^s(p_{1,\mu_0})$ and $W^u(p_{2,\mu_0})$ have a quadratic tangency $q_{\mu_0}$ at $\mu_0\in J$.
We say that the tangency $q_{\mu_0}$ {\it unfolds generically} with respect to $\{\psi_\mu\}_{\mu\in J}$ if there exist local 
coordinates $(x,y)$ on $\mathcal{N}_\mu$ and $C^2$ functions $\alpha_\mu(x)$ which $C^2$ depend on $\mu$ and satisfy the following conditions, where $\{\mathcal{N}_\mu\}$ is a $C^2$ family of small open neighborhoods of $q_\mu$ in $\mathbb{R}^2$.
\begin{itemize}
\item
$\alpha_{\mu_0}(x)$ satisfies the condition $(\ref{eqn_tangency})_{n=1}$ and $\alpha_{\mu_0}(0)=q_{\mu_0}$.
\item
$\{(x, y) ; y=0\}\subset W^{s}(p_{1,\mu})$ and $\{(x, y) ;  y=\alpha_\mu(x)\}\subset W^{u}(p_{2,\mu})$ for any $\mu\in J$ near $\mu_0$. 
\item
For the two variable function $\alpha(\mu,x):=\alpha_\mu(x)$,
\begin{equation}\label{qt_unfold}
\frac{\partial \alpha}{\partial \mu}(\mu_0,0)\neq 0.
\end{equation}
\end{itemize}
It is not hard to see that the definition of this generic condition is independent of the choice of the coordinate neighborhoods $\mathcal{N}_\mu$ satisfying the conditions as above.

\begin{figure}[hbt]
\begin{center}
\scalebox{0.8}{\includegraphics[clip]{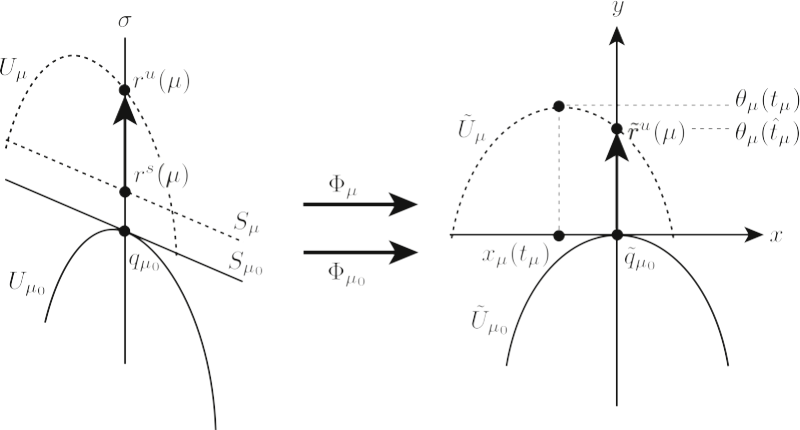}}
\caption{}
\label{fig_11}
\end{center}
\end{figure}

Now, we will present a practical condition equivalent to (\ref{qt_unfold}) which works under general coordinates on a neighborhood of $q_{\mu_0}$.
Suppose that $\tilde {\mathcal N}_\mu$ are coordinates with respect to which there exists a continuation $S_\mu$ of curves in $W^s(p_{1,\mu})$ with $q_{\mu_0}\in \mathrm{Int}S_{\mu_0}$ which are represented by the graphs of $C^2$ functions $\eta_\mu(x)$ of $x$ with $|x|\leq \delta$ for some $\delta>0$, that is,
$$S_\mu=\{(x,\eta_\mu(x));\, |x|\leq \delta\}.$$
We call $\eta_\mu$ the \emph{graph function} of $S_\mu$.
Let $U_\mu$ be a continuation of curves in $W^u(p_{2,\mu})$ with $q_{\mu_0}\in \mathrm{Int}U_{\mu_0}$, and $\sigma$ a vertical segment passing through $S_{\mu_0}$ at $q_{\mu_0}$, see Figure \ref{fig_11}.
%%%%%%%%%%%%%%%%%%%%%%%%%
%%%%%%%%%%%%%%%%%%%%%%
The intersection $S_\mu\cap \sigma$ (resp.\ $U_\mu\cap \sigma$) defines a continuation of points $r^s(\mu)$ (resp.\ $r^u(\mu)$).
We denote the velocity vectors of $r^s(\mu)$ and $r^u(\mu)$ at $\mu=\mu_0$ by $\boldsymbol{v}_\mu^{s,\perp}(q_{\mu_0})$, $\boldsymbol{v}_\mu^{u,\perp}(q_{\mu_0})$ respectively.
That is,
$$\boldsymbol{v}_\mu^{s,\perp}(q_{\mu_0})=\frac{d}{d\mu}r^s(\mu_0),\quad \boldsymbol{v}_\mu^{u,\perp}(q_{\mu_0})=\frac{d}{d\mu}r^u(\mu_0).$$

Let $\Phi_\mu$ be the coordinate change of $\tilde {\mathcal N}_\mu$ defined by $\Phi_\mu(x,y)=(x,y-\eta_\mu(x))$.
Then $\tilde S_\mu=\Phi_\mu(S_\mu)$ is contained in the $x$-axis.
Let $(x_\mu(t),y_\mu(t))$  be a $C^2$ regular curve parametrization of $U_\mu$ which $C^2$ depends on $\mu$ and 
such that the curve passes through $q_{\mu_0}$ at $t=0$, 
where a plane curve $c(t)$ being \emph{regular} means that $dc(t)/dt\neq (0,0)$ for any $t$.
Set
\begin{equation}\label{theta_mu}
\theta_\mu(t)=y_\mu(t)-\eta_\mu(x_\mu(t)).
\end{equation}
Then $(x_\mu(t),\theta_\mu(t))$ is a parametrization of $\tilde U_\mu=\Phi_\mu(U_\mu)$.
Since $\tilde U_\mu$ is a quadratic curve for any $\mu$ close to $\mu_0$, there exists a unique $t_\mu$ near $0$ such that $\theta_\mu(t)$ has an extremal point at $t=t_\mu$ which $C^1$ depends on $\mu$.
Similarly, since $\tilde U_\mu$ meets the $y$-axis transversely in a single point $\tilde r^u(\mu)$, there exists a unique $\hat t_\mu$ near $0$ with $(0,\theta_\mu(\hat t_\mu))=\tilde r^u(\mu)$ which $C^2$ depends on $\mu$.
The generic condition (\ref{qt_unfold}) for $\tilde U_{\mu_0}$ with respect to the new coordinate is ${d \theta_\mu(\hat t_\mu)}/{d \mu}(\mu_0)\neq 0$.
Since $\tilde r^u(\mu)=r^u(\mu)-r^s(\mu)$ as a vector,
$$\Bigl(0,\frac{d\theta_\mu(\hat t_\mu)}{d\mu}(\mu_0)\Bigr)=\frac{d\tilde r^u}{d\mu}(\mu_0)=
\boldsymbol{v}_\mu^{u,\perp}(q_{\mu_0})-\boldsymbol{v}_\mu^{s,\perp}(q_{\mu_0}).$$
From the definitions as above, $t_{\mu_0}=\hat t_{\mu_0}$.
Thus, $t_{\mu}-\hat t_{\mu}=O(\Delta \mu)$ for $\mu=\mu_0+\Delta\mu$.
Since $\theta_\mu(t)$ has an extremal value at $t=t_\mu$,
$$\theta_\mu(\hat t_\mu)-\theta_\mu(t_\mu)=O((\hat t_\mu-t_\mu)^2)=O(\Delta\mu^2).$$
This shows that $d\theta_\mu(\hat t_\mu)/d\mu|_{\mu=\mu_0}=d\theta_\mu(t_\mu)/d\mu|_{\mu=\mu_0}$.
Hence the generic condition (\ref{qt_unfold}) is equivalent to
\begin{equation}\label{qt_unfold_3}
\boldsymbol{v}_\mu^{u,\perp}(q_{\mu_0})-\boldsymbol{v}_\mu^{s,\perp}(q_{\mu_0})=\Bigl(0,\frac{d\theta_\mu(t_\mu)}{d\mu}(\mu_0)\Bigr)\neq (0,0).
\end{equation}

\subsection{Generically unfolding cubic tangencies}\label{S_gct}

Suppose that $\psi$ is a $C^3$ diffeomorphism of $\mathbb{R}^2$ with a saddle fixed point $p$.
A cubic homoclinic tangency $q$ of $\psi$ associated with $p$ is said to \textit{unfold generically} with respect to a two-parameter family $\{\psi_{u,v}\}$ in $\mathrm{Diff}^3(\mathbb{R}^2)$ with $\psi_{0,0}=\psi$ if there exist $(u,v)$-dependent local coordinates $(x, y)$ on a neighborhood of $q$ with $q = (0,0)$ such that $W^{s}(p_{u,v}) = \{(x, y) ; y=0\}$ and 
$W^{u}(p_{u,v}) = \{(x, y) ; y=y_{u,v}(x)\}$ and $y_{u,v}(x)=y(u,v,x)$ is a $C^3$ function satisfying
\begin{equation}\label{generic condition}
(\partial_{u}y \cdot \partial_{v x}y-\partial_{v}y\cdot  \partial_{u x}y)(0,0,0)\neq 0,
\end{equation}
where $\{p_{u,v}\}$ is a continuation of saddle fixed points of $\psi_{u,v}$ with $p_{0,0}=p$. 
Since $y_{0,0}(0)=y_{0,0}'(0)=y_{0,0}''(0)=0$, $y_{\mu,\nu}$ has the Taylor expansion
\begin{equation}\label{Taylor_exp}
y_{u,v}(x)=a_1u+a_2v+a_3ux+a_4vx+a_5uv+h(x,u,v),
\end{equation}
where $a_1,\dots,a_5$ are constants and $h(x,u,v)$ is a $C^3$ function with
\begin{equation}\label{condition_h}
h=\partial_u h=\partial_v h=\partial_x h=\partial_{ux}h=\partial_{vx} h=\partial_{xx}h=0
\end{equation}
at $(x,u,v)=(0,0,0)$.
Then the generic condition (\ref{generic condition}) is rewritten as follows.
\begin{equation}\label{gc_2}
a_1a_4-a_2a_3\neq 0.
\end{equation}
Let $F:(u,v)\mapsto (\hat u,\hat v)$ is a $C^3$-diffeomorphism with $F(0,0)=(0,0)$, and let
$$
y_{F^{-1}(\hat u,\hat v)}(x)=b_1\hat u+b_2\hat v+b_3\hat ux+b_4\hat vx+b_5\hat u\hat v+\hat h(x,F^{-1}(\hat u,\hat v))
$$
be the expansion of $y_{F^{-1}(\hat u,\hat v)}$, where $\hat h$ is a $C^3$ function satisfying 
the condition as (\ref{condition_h}).
Then we have
$$\begin{pmatrix}a_1&a_2\\a_3&a_4\end{pmatrix}=\begin{pmatrix}b_1&b_2\\b_3&b_4\end{pmatrix}DF_{(0,0)}.$$
This equation implies the following.

\begin{lemma}\label{l_ct_00}
With the notation as above, $q=(0,0)$ is a cubic tangency unfolding generically with respect to $\{\psi_{u,v}\}$ if and only if it is one unfolding generically with respect to $\{\psi_{F^{-1}(\hat u,\hat v)}\}$.
\end{lemma}

Now, we show that this generic condition is preserved under coordinate changes of the $xy$-plane fixing the $x$-axis as a set.
Let $U$ be a small neighborhood of $(0,0)$ in the $uv$-space.
Suppose that $\{\Phi_{u,v}\}_{(u,v)\in U}$ is a two-parameter family of $C^3$ diffeomorphisms of the $xy$-plane which $C^3$ depends on $(u,v)$ and such that each $\Phi_{u,v}$ fixes the $x$-axis as a set and $\Phi_{0,0}(0,0)=(0,0)$.
Let $\rho_{u,v}$ is a continuation of curves in $W^u(p_{u,v})$ with $\mathrm{Int}\rho_{u,v}\ni q$.
We set
$$\tilde \rho_{u,v}=\Phi_{u,v}(\rho_{u,v}),\quad \tilde \psi_{u,v}=\Phi_{u,v}\circ \psi_{u,v}\circ \Phi_{u,v}^{-1}.$$

\begin{lemma}\label{l_ct_0}
With the notation as above, if $q=(0,0)$ is a cubic tangency of $\rho_{0,0}$ and the $x$-axis which unfolds generically with respect to $\{\psi_{u,v}\}$, then $\tilde q=\Phi_{0,0}(q)$ is also a cubic tangency of $\tilde \rho_{0,0}$ and the $x$-axis which unfolds generically with respect to $\{\tilde \psi_{u,v}\}$.
\end{lemma}
\begin{proof}
Since $\Phi_{u,v}$ preserves the $x$-axis, $\Phi_{u,v}$ is represented as
$$\Phi_{u,v}(x,y)=(\beta_{u,v}(x,y),y\gamma_{u,v}(x,y)),$$
where $\beta_{u,v}$ (resp.\ $\gamma_{u,v}$) is a $C^3$ (resp.\ $C^2$) function.
Moreover, the condition $\Phi_{0,0}(0,0)=(0,0)$ implies 
$$\beta_{0,0}(0,0)=0.$$
Since the differential of $\Phi_{u,v}$ is
$$D\Phi_{u,v}=\begin{pmatrix}
\partial_x \beta_{u,v}&\partial_y\beta_{u,v}\\
y\partial_x\gamma_{u,v}&\gamma_{u,v}+y\partial_y \gamma_{u,v}
\end{pmatrix},
$$
$\det (D\Phi_{0,0}(0,0))=\partial_x\beta_{0,0}(0,0)\cdot \gamma_{0,0}(0,0)$.
Since $\Phi_{0,0}$ is a diffeomorphism,
$$b=\partial_x \beta_{0,0}(0,0)\neq 0,\quad c=\gamma_{0,0}(0,0)\neq 0.$$
The curve $\tilde \rho_{u,v}$ is parametrized as
$$\Phi_{u,v}(x,y_{u,v}(x))=(\beta_{u,v}(x,y_{u,v}(x)),y_{u,v}(x)\cdot \gamma_{u,v}(x,y_{u,v}(x))).$$
Set $\tilde x=\tilde x(x,u,v)=\beta_{u,v}(x,y_{u,v}(x))$.
Then $\tilde x(0,0,0)=\beta_{0,0}(0,0)=0$.
Differentiating $\tilde x$ by $x$,
$$\partial_x \tilde x=\partial_x \beta_{u,v}(x,y_{u,v}(x))+\partial_y \beta_{u,v}(x,y_{u,v}(x))\partial_x y_{u,v}(x).$$
Since $\partial_x y_{0,0}(0)=0$, $\partial_x \tilde x(0,0,0)=b\neq 0$.
Thus, $\tilde x(x,u,v)$ has a local inverse function $x=\delta_{u,v}(\tilde x)$ with $\delta_{0,0}(0)=0$ defined for any $(u,v)$ near $(0,0)$ and any $\tilde x$ near $0$.
It follows that
$$\Phi_{u,v}(x,y_{u,v}(x))=(\tilde x,\tilde y_{u,v}(\tilde x))=(\tilde x,\tilde w_{u,v}(\tilde x)\cdot \tilde \gamma_{u,v}(\tilde x)),$$
where $\tilde w_{u,v}(\tilde x)=y_{u,v}\circ \delta_{u,v}(\tilde x)$ and $\tilde \gamma_{u,v}(\tilde x)=\gamma_{u,v}(\delta_{u,v}(\tilde x),\tilde w_{u,v}(\tilde x))$.
Then
\begin{equation}\label{partial_u_y}
\partial_u \tilde w_{u,v}(\tilde x)=(\partial_u y_{u,v})(\delta_{u,v}(\tilde x))+(\partial_x y_{u,v})(\delta_{u,v}(\tilde x))\partial_u \delta_{u,v}(\tilde x).
\end{equation}
Since $\partial_x y_{0,0}(0)=0$, $\partial_u \tilde w_{u,0}(0)|_{u=0}=a_1$.
We have as well $\partial_v \tilde w_{0,v}(0)|_{v=0}=a_2$.
From $\tilde y_{u,v}(\tilde x)=\tilde w_{u,v}(\tilde x)\cdot \tilde \gamma_{u,v}(\tilde x)$,
\begin{equation}\label{y_cdot_gamma}
\partial_u (\tilde y_{u,v}(\tilde x))=\partial_u \tilde w_{u,v}(\tilde x)\cdot \tilde\gamma_{u,v}(\tilde x)+\tilde w_{u,v}(\tilde x)\cdot \partial_u \tilde\gamma_{u,v}(\tilde x).
\end{equation}
Since $\tilde \gamma_{0,0}(0)=\gamma_{0,0}(0)=c$ and $\tilde w_{0,0}(0)=0$, $\partial_u (\tilde y_{u,0}(0))|_{u=0}=a_1c.$
A similar argument shows $\partial_v (\tilde y_{0,v}(0))|_{v=0}=a_2c.$

Differentiating the both sides of (\ref{partial_u_y}) by $\tilde x$ and putting $(u,\tilde x)=(0,0)$, we have
$$\partial_{u\tilde x}(\tilde w_{u,0}(\tilde x))|_{(u,\tilde x)=(0,0)}=\frac{a_3}{b}.$$
Then, from this equation together with the differentiation of (\ref{y_cdot_gamma}) by $\tilde x$, we have
$$
\partial_{u\tilde x} (\tilde y_{u,0}(\tilde x))|_{(u,\tilde x)=(0,0)}=\frac{a_3c}{b}+a_1d,
$$
where $d=\partial_{\tilde x}\tilde \gamma_{0,0}(\tilde x)|_{\tilde x=0}$.
Similarly,
$$
\partial_{v\tilde x} (\tilde y_{0,v}(\tilde x))|_{(v,\tilde x)=(0,0)}=\frac{a_4c}{b}+a_2d.
$$
By using the equalities as above,
\begin{align*}
&\partial_u (\tilde y_{u,0}(0))\cdot 
\partial_{v\tilde x} (\tilde y_{0,v}(\tilde x))|_{(u,v,\tilde x)=(0,0,0)}-
\partial_v (\tilde y_{0,v}(0))\cdot \partial_{u\tilde x} (\tilde y_{u,0}(\tilde x))|_{(u,v,\tilde x)=(0,0,0)}\\
&\qquad\qquad\qquad=a_1c\Bigl(\frac{c}{b}a_4+a_2d\Bigr)-a_2c\Bigl(\frac{c}{b}a_3+a_1d\Bigr)\\
&\qquad\qquad\qquad=\frac{c^2}{b}(a_1a_4-a_2a_3)\neq 0.
\end{align*}
It follows that $(0,0)$ is a cubic tangency of $\tilde \rho_{0,0}$ and the $\tilde x$-axis which unfolds generically with respect to $\{\tilde \psi_{u,v}\}$.
\end{proof}

\section{Gonchenko-Shilnikov-Turaev's approach to cubic tangencies; revisited}\label{S_GST}

In series of papers \cite{GST1}\,--\,\cite{GST7}, Gonchenko, Shilnikov and Turaev studied extensively generic 
homoclinic tangencies of higher order for two-dimensional diffeomorphisms.
In this section, we present some of their results on cubic homoclinic tangencies in a form suitable to 
our proofs.

\subsection{Existence of generic cubic homoclinic tangencies}\label{SS_exist_cubic}

The following lemma is a special case of \cite[Lemma 5]{GST7}.

\begin{lemma}\label{G_Lemma_5}
Suppose that $\psi$ is a $C^r$-diffeomorphism with $4\leq r\leq \omega$ on the plane $\mathbb{R}^2$ with two saddle fixed points $p^+,p^-$.
Let $\{\psi_{\mu,\nu}\}$ be a two-parameter family of $C^r$ diffeomorphisms on $\mathbb{R}^2$ with $\psi_{0,0}=\psi$ and $C^{r-2}$ depending on $\mu,\nu$.
Let $p^\pm_{\mu,\nu}$ be continuations of saddle fixed points of $\psi_{\mu,\nu}$ with $p_{0,0}^+=p^+$, $p_{0,0}^-=p^-$.
Suppose the following conditions.
\begin{enumerate}[\rm (i)]
\item
$W^s(p^+)$ and $W^u(p^-)$ have a quadratic heteroclinic tangency $q^+$ unfolding generically with respect to $\{\psi_{\mu,0}\}$.
\item
$W^s(p^-)$ and $W^u(p^+)$ have a quadratic heteroclinic tangency $q^-$ unfolding generically with respect to $\{\psi_{0,\nu}\}$.
\end{enumerate}
Then there exists an element $(\mu^*,\nu^*)\neq (0,0)$ in the $\mu\nu$-space arbitrarily close to $(0,0)$ such that $\psi_{\mu^*,\nu^*}$ has a cubic homoclinic tangency $r_{\mu^*,\nu^*}$ associated with $p^+_{\mu^*,\nu^*}$ which unfolds generically with respect to $\{\psi_{\mu,\nu}\}$.
\end{lemma}

\begin{figure}[hbt]
\begin{center}
\scalebox{0.8}{\includegraphics[clip]{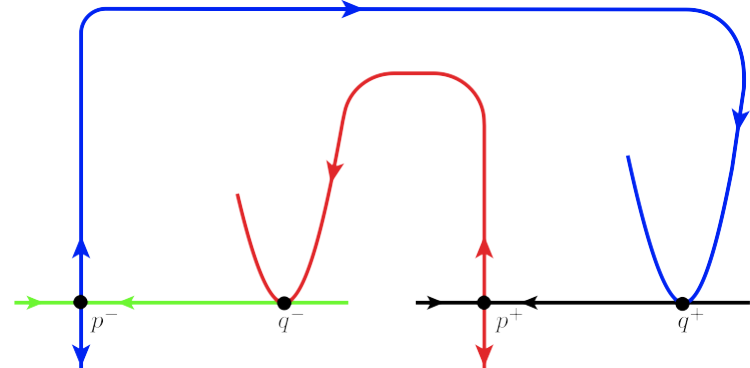}}
\caption{}
\label{fig_21}
\end{center}
\end{figure}

The situation of Lemma \ref{G_Lemma_5} is illustrated in Figure \ref{fig_21}.
%%%%%%%%%%%%%%%%%%%%%%%%%
%%%%%%%%%%%%%%%%%%%%%%
Gonchenko et al \cite{GST7} show that the existence of a cubic homoclinic tangency under the assumptions of this lemma.
But they did not present the fact that the cubic tangency $r_{\mu^*,\nu^*}$ unfolds generically.
Here we will review the proof of Lemma \ref{G_Lemma_5} and show the unfoldingness, which is crucial in 
the proof of Theorem \ref{main_a}.
Note that our notations $q^-,q^+,\mu,\nu,p^+,p^-$ here correspond respectively to the notations $M_+,M_-,\nu,\mu_0,O_1=O_3,O_2$ in 
\cite[Lemma 5]{GST7}.

\begin{figure}[hbt]
\begin{center}
\scalebox{0.8}{\includegraphics[clip]{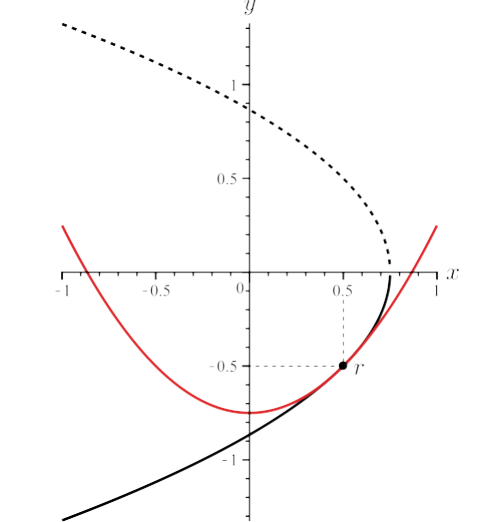}}
\caption{The curve $y=-\frac34+x^2$ and the minus branch $y=-\sqrt{\frac34-x}$ of the curve $x=\frac34-y^2$ have a cubic tangency at $\bigl(\frac12,-\frac12\bigr)$.}
\label{fig_22}
\end{center}
\end{figure}

Figure \ref{fig_22} illustrates a typical example of a cubic tangency of two quadratic curves.
%%%%%%%%%%%%%%%%%%%%%%%%%
A basic idea of the proof of Lemma \ref{G_Lemma_5} is to reparametrize $x,y,\mu,\nu$ so as to approximate $W^u(p_{\mu,\nu}^+)$ and $W^s(p_{\mu,\nu}^+)$ near $p_{\mu,\nu}^-$ by the curves $y=-\frac34+x^2$ and $x=\frac34-y^2$.

\begin{proof}[Proof of Lemma \ref{G_Lemma_5}]
Replacing $q^+$ by $\psi^{-n}(q^+)$ and $q^-$ by $\psi^{n}(q^-)$ for some $n\in \mathbb{N}$, we may assume that 
both $q^+$, $q^-$ are contained in an arbitrarily small neighborhood $U$ of $p^-$ in $\mathbb{R}^2$.
Let $\lambda=\lambda_{\mu,\nu}, \gamma=\gamma_{\mu,\nu}$ be the eigenvalues of $D\psi_{\mu,\nu}(p_{\mu,\nu}^-)$ with $0<|\lambda|<1<|\gamma|$.
By Afraimovich-Shilnikov \cite{AS} together with Gonchenko-Shilnikov \cite{GS, GS2}, one can choose $C^{r-1}$ coordinate systems $(x,y)$ on $U$ and reparametrize $(\mu,\nu)$ if necessary so as to satisfy the following conditions (i)--(iv), where the equations (\ref{G_3.23}), (\ref{G_3.24}) correspond respectively to (3.23), (3.24) in \cite{GST7}.
\begin{enumerate}[(i)]
\item
The $x$-axis in $U$ is $W^s_{\mathrm{loc}}(p_{\mu,\nu}^-)$ and the $y$-axis is $W^u_{\mathrm{loc}}(p_{\mu,\nu}^-)$ and $q^-=(x_0,0)$, $q^+=(0,y_0)$ for some 
non-zero constants $x_0,y_0$.
In particular, $p_{\mu,\nu}^-$ is the origin with respect to the coordinate. 
\item
$W^u(p^+_{\mu,\nu})$ contains a curve $C^{r-1}$ depending on $(\mu,\nu)$ and defined by the equation
\begin{equation}\label{G_3.23}
y=\nu+d(x-x_0)^2+o\bigl((x-x_0)^2\bigr)
\end{equation}
for some non-zero constant $d$.
In particular, the curves for $\nu=0$ contain $q^-$.
\item
$W^s(p^+_{\mu,\nu})$ contains a curve $C^{r-1}$ depending on $(\mu,\nu)$ and defined by the equation
\begin{equation}\label{G_3.24}
x=\mu+\hat d(y-y_0)^2+o\bigl((y-y_0)^2\bigr)
\end{equation}
for some non-zero constant $\hat d$.
In particular, the curves for $\mu=0$ contain $q^+$.
\item
Let $k$ be any positive integer.
If $(x^{(i)},y^{(i)}):=\psi_{\mu,\nu}^i(x,y)$ are contained in $U$ for all $i=0,1,\dots,k$, then 
$$x^{(k)}=\lambda^kx^{(0)}+\lambda^k\xi_k(x^{(0)},y^{(k)}),\quad 
y^{(0)}=\gamma^{-k}y^{(k)}+\gamma^{-k}\eta_k(x^{(0)},y^{(k)}),$$ 
where $\xi_k$, $\eta_k$ are functions satisfying 
\begin{equation}\label{G_2.7}
\Vert \xi_k,\eta_k\Vert_{C^{r-1},C^{r-2}}=o(1)_{k\rightarrow \infty}.
\end{equation}
\end{enumerate}
Here (\ref{G_2.7}) means that these functions themselves and the derivatives of them up to order $r-1$ with 
respect to the variables and up to order $r-2$ with respect to the parameters converge uniformly to 
zero as $k\rightarrow \infty$.

By the coordinate change $(x,y)\mapsto (X,Y)$ in a small neighborhood of $(x_0,0)$ defined by $X=x-x_0$, $Y=y^{(k)}-y_0$, the curve (\ref{G_3.23}) and the $\psi_{\mu,\nu}^{-k}$-image of the 
curve (\ref{G_3.24}) 
are given by
\begin{align}
&\gamma^{-k}(Y+y_0+\eta_k(X+x_0,Y+y_0))=\nu+dX^2+o(X^2),
\tag{\ref{G_3.23}a}\label{G_3.25}\\
&\lambda^k(X+x_0+\xi_k(X+x_0,Y+y_0))=\mu+\hat d Y^2+o(Y^2).\tag{\ref{G_3.24}a}\label{G_3.26}
\end{align}
%%%%%%%%%%%%%%%%%%%%%%%%%
\begin{figure}[hbt]
\begin{center}
\scalebox{0.8}{\includegraphics[clip]{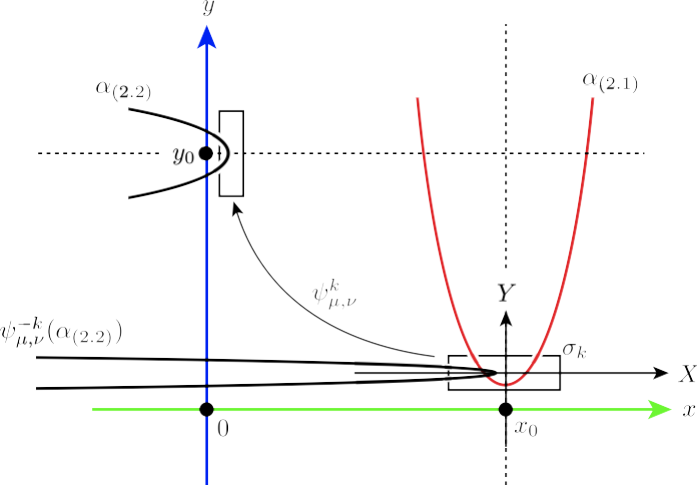}}
\caption{
$\sigma_k$ is a strip near $(x_0,0)$ as illustrated in the figure such that $\psi_{\mu,\nu}^k(\sigma_k)$ is a strip near $(0,y_0)$.
The $X$-axis is contained in the $\psi_{\mu,\nu}^{-k}$-image of the horizontal line $y=y_0$ and 
the $Y$-axis is in the vertical line $x=x_0$.
$\alpha_{(2,1)}$, $\alpha_{(2,2)}$ are curves represented by the equations (2.1), (2.2), respectively.
}
\label{fig_21_1}
\end{center}
\end{figure}
%%%%%%%%%%%%%%%%%%%%%%
See Figure \ref{fig_21_1}.
By rearranging (\ref{G_3.25}), (\ref{G_3.26}) and applying (\ref{G_2.7}), we have
\begin{align}
\gamma^{-k}Y&=\nu-\gamma^{-k}y_0+o(\gamma^{-k})+dX^2+o(X^2),\tag{\ref{G_3.23}b}\label{G_3.29}\\
\lambda^k X&=\mu-\lambda^kx_0+o(\lambda^k)+\hat d Y^2+o(Y^2).\tag{\ref{G_3.24}b}\label{G_3.30}
\end{align}
Next we consider the $k$-dependent reparametrization $(X,Y,\mu,\nu)\mapsto (\tilde x,\tilde y,\tilde \mu,\tilde \nu)$ defined as
\begin{align*}
Y&=(\hat d^2 d)^{-1/3}\lambda^{2k/3}\gamma^{-k/3}\tilde y,\quad
X=-(\hat d d^2)^{-1/3}\lambda^{k/3}\gamma^{-2k/3}\tilde x,\\
\tilde\mu&= (\mu-\lambda^kx_0+o(\lambda^k))d^{2/3}\hat d^{1/3}\lambda^{-4k/3}\gamma^{2k/3},\\
\tilde\nu&= (\nu-\gamma^{-k}y_0+o(\gamma^{-k}))d^{1/3}\hat d^{2/3}\lambda^{-2k/3}\gamma^{4k/3}.
\end{align*}
Then the curve (\ref{G_3.23}) and the $\psi_{\mu,\nu}^{-k}$-image of the curve (\ref{G_3.24}) are represented respectively by 
\begin{align}
\tilde y&=\tilde\nu+\tilde x^2+o(1)_{k\rightarrow \infty},\tag{\ref{G_3.23}c}\label{G_3.32}\\
\tilde x&=\tilde\mu-\tilde y^2+o(1)_{k\rightarrow \infty},\tag{\ref{G_3.24}c}\label{G_3.33}
\end{align}
where the term $o(1)_{k\rightarrow \infty}$ in (\ref{G_3.32}) is $O(X)=O(\lambda^{k/3}\gamma^{-k/3})$ and
the term $o(1)_{k\rightarrow \infty}$ in (\ref{G_3.33}) is $O(Y)=O(\lambda^{k/3}\gamma^{-k/3})$.  
Note that, if $(\tilde \mu, \tilde \nu)$ ranges a bounded domain, then 
$(\mu,\nu)$ converges to $(0,0)$ as $k\rightarrow \infty$.

The curve (\ref{G_3.32}) and the minus branch 
$\tilde y=-\sqrt{\tilde\mu-\tilde x+o(1)_{k\rightarrow \infty}}$ of (\ref{G_3.33}) have 
a cubic tangency if and only if the curve
\begin{equation}\label{minus_branch}
\tilde y(\tilde \mu,\tilde\nu,\tilde x)=\tilde\nu+\tilde x^2+\sqrt{\tilde\mu-\tilde x+o(1)_{k\rightarrow \infty}}
+o(1)_{k\rightarrow \infty}
\end{equation}
and the $\tilde x$-axis have a cubic tangency.
In fact, Gonchenko et al proved that, for all sufficiently large $k\in \mathbb{N}$, there exist 
$\tilde x^*=\frac12+o(1)_{k\rightarrow \infty}$, $\tilde\nu^*=-\frac34+o(1)_{k\rightarrow \infty}$, $\tilde\mu^*=\frac34+o(1)_{k\rightarrow \infty}$ such that 
the curve (\ref{minus_branch})$_{(\tilde\mu,\tilde\nu)=(\tilde \mu^*,\tilde \nu^*)}$ and the $\tilde x$-axis 
have a cubic tangency at $(\tilde x^*,0)$ and hence the curves (\ref{G_3.32})$_{(\tilde\mu,\tilde\nu)=(\tilde \mu^*,\tilde \nu^*)}$ and (\ref{G_3.33})$_{(\tilde\mu,\tilde\nu)=(\tilde \mu^*,\tilde \nu^*)}$ also have 
a cubic tangency $\tilde r_{\tilde \mu^*,\tilde\nu^*}$ at $\bigl(\frac12+o(1)_{k\rightarrow \infty},-\frac12+o(1)_{k\rightarrow \infty}\bigr)$.
Then the point $r_{\mu^*,\nu^*}$ in the $xy$-plane corresponding to $\tilde r_{\tilde \mu^*,\tilde\nu^*}$ in the $\tilde x\tilde y$-plane is a cubic tangency of the curve (\ref{G_3.23})$_{(\mu,\nu)=(\mu^*,\nu^*)}$ and the $\psi_{\mu^*,\nu^*}^{-k}$-image of (\ref{G_3.24})$_{(\mu,\nu)=(\mu^*,\nu^*)}$, 
see Figure \ref{fig_23}.
%%%%%%%%%%%%%%%%%%%%%%%%%
\begin{figure}[hbt]
\begin{center}
\scalebox{0.8}{\includegraphics[clip]{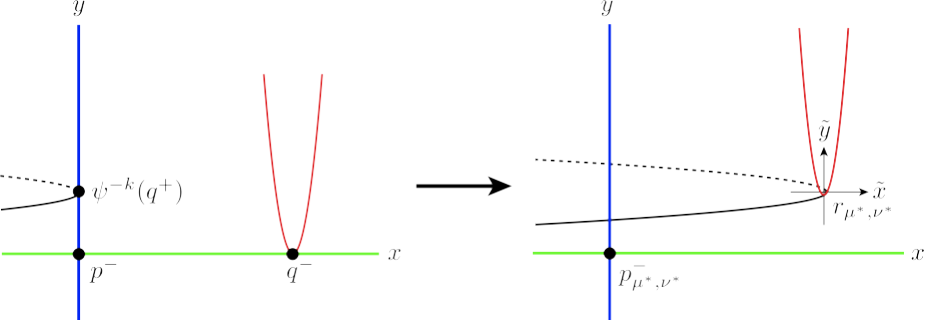}}
\caption{A small neighborhood of $r_{\mu^*,\nu^*}$ in the right hand side has a situation similar 
up to scale to that in Figure \ref{fig_22}.}
\label{fig_23}
\end{center}
\end{figure}
%%%%%%%%%%%%%%%%%%%%%%
Since $|\tilde \mu^*|$, $|\tilde \nu^*|<1$ for all sufficiently large $k$, one can suppose that the corresponding parameter value 
$(\mu^*,\nu^*)$ is arbitrarily close to $(0,0)$.

By (\ref{minus_branch}), we have
\begin{align*}
\partial_{\tilde\mu} \tilde y(\tilde\mu^*,\tilde\nu^*,\tilde x^*)&=
\frac12(\tilde\mu^*-\tilde x^*)^{-\frac12}+o(1)_{k\rightarrow \infty}=1+o(1)_{k\rightarrow \infty},\\
\partial_{\tilde\mu\tilde x} \tilde y(\tilde\mu^*,\tilde\nu^*,\tilde x^*)&=
\frac14(\tilde\mu^*-\tilde x^*)^{-\frac32}+o(1)_{k\rightarrow \infty}=2+o(1)_{k\rightarrow \infty},\\
\partial_{\tilde\nu} \tilde y(\tilde\mu^*,\tilde\nu^*,\tilde x^*)&=1+o(1)_{k\rightarrow \infty},\  
\partial_{\tilde\nu\tilde x} \tilde y(\tilde\mu^*,\tilde\nu^*,\tilde x^*)=o(1)_{k\rightarrow \infty}.
\end{align*}
It follows that
$$(\partial_{\tilde\mu}\tilde y\cdot \partial_{\tilde\nu\tilde x}\tilde y-\partial_{\tilde\nu}\tilde y\cdot \partial_{\tilde\mu\tilde x}\tilde y)(\tilde\mu^*,\tilde\nu^*,\tilde x^*)=-2+o(1)_{k\rightarrow \infty}\neq 0.$$
Thus the cubic tangency $\tilde r_{\tilde\mu^*,\tilde\nu^*}$ unfolds generically.
By Lemmas \ref{l_ct_00} and \ref{l_ct_0}, the cubic tangency $r_{\mu^*,\nu^*}$ also unfolds
generically.
\end{proof}

\begin{figure}[hbt]
\begin{center}
\scalebox{0.8}{\includegraphics[clip]{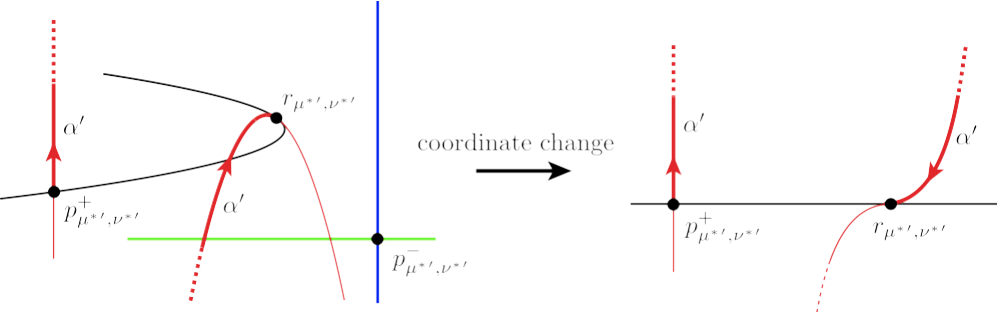}}
\caption{A cubic tangency of type I.}
\label{fig_25}
\end{center}
\end{figure}

\begin{definition}\label{dfn_cub_type}
Let $r$ be a cubic homoclinic tangency of $\varphi \in \mathrm{Diff}^3(\mathbb{R}^2)$ associated with 
a saddle fixed point $p$ of $\varphi$ and $\alpha$ the arc in $W^u(p)$ connecting $p$ with $r$.
We say that the tangency is of \emph{type I} if small neighborhoods of $p$ and $r$ in $\alpha$ are contained in the same side of $W^s(p)$, see the right hand side of Figure \ref{fig_25}.
Otherwise it is called of \emph{type II}, see Figure \ref{fig_24}.
\end{definition}

\begin{lemma}\label{type_1}
With the notation as in Lemma \ref{G_Lemma_5}, suppose moreover that the unstable eigenvalue of $D\psi(p^+)$ is positive and the both eigenvalues of $D\psi(p^-)$ are negative.
Then, the cubic homoclinic tangency $r_{\mu^*,\nu^*}$ given in Lemma \ref{G_Lemma_5} can be chosen so 
that it is of type I.
\end{lemma}
\begin{proof}
Let $\alpha_0$ be the arc in $W^u(p^+)$ connecting $p^+$ with $q^-$.
Recall that $r_{\mu^*,\nu^*}$ is obtained by perturbing the curve (\ref{G_3.23})$_{(\mu,\nu)=(0,0)}$ 
and the $\psi^{-k}$-image of (\ref{G_3.24})$_{(\mu,\nu)=(0,0)}$.
Suppose that $r_{\mu^*,\nu^*}$ is of type II.
The arc $\alpha$ in $W^u(p_{\mu^*,\nu^*}^+)$ connecting $p_{\mu^*,\nu^*}^+$ with $r_{\mu^*,\nu^*}$ 
is obtained by slightly deforming $\alpha_0$ in $\mathbb{R}^2$, see Figure \ref{fig_24}.
%%%%%%%%%%%%%%%%%%%%%%%%%
\begin{figure}[hbt]
\begin{center}
\scalebox{0.8}{\includegraphics[clip]{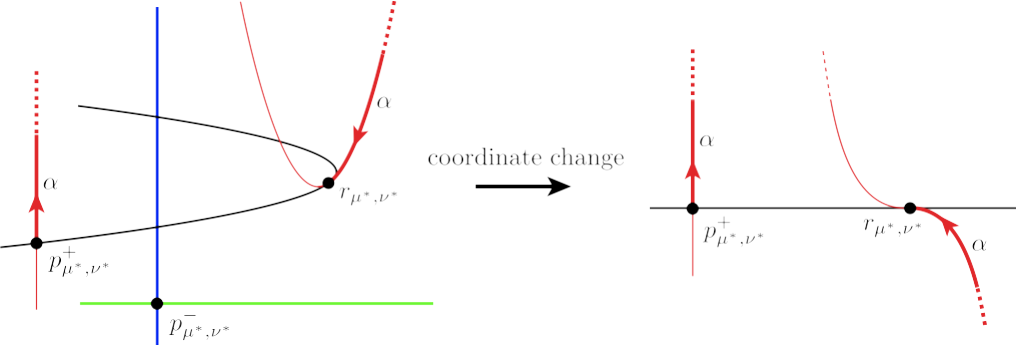}}
\caption{A cubic tangency of type II.}
\label{fig_24}
\end{center}
\end{figure}
%%%%%%%%%%%%%%%%%%%%%%
The arc $\alpha_0'=\psi(\alpha_0)$ connects $p^+$ with $\psi(q^-)$ in $W^u(p^+)$.
Since the unstable eigenvalue of $D\psi(p^+)$ is positive, $\alpha_0'$ contains $\alpha_0$.
Since the both eigenvalues of $D\psi(p^-)$ are negative, the curve (\ref{G_3.23})$_{(\mu,\nu)=(0,0)}$ and 
its $\psi$-image lie in mutually opposite quadrants with respect to the coordinate $(x,y)$ used in 
Lemma \ref{G_Lemma_5}.
Again by applying Lemma \ref{G_Lemma_5}, one can show that there exist small parameter values $\mu^*{}'$, $\nu^*{}'$ with $\mu^*{}'\mu^*<0$, $\nu^*{}'\nu^*<0$ such that the $\psi_{\mu^*{}',\nu^*{}'}$-image of (\ref{G_3.23})$_{(\mu,\nu)=({\mu^*{}',\nu^*{}'})}$ 
and the $\psi_{\mu^*{}',\nu^*{}'}^{-k}$-image of (\ref{G_3.24})$_{(\mu,\nu)=({\mu^*{}',\nu^*{}'})}$ have a cubic tangency $r_{\mu^*{}',\nu^*{}'}$.
The arc $\alpha'$ in $W^u(p_{\mu^*{}',\nu^*{}'}^+)$ connecting $p_{\mu^*{}',\nu^*{}'}^+$ with $r_{\mu^*{}',\nu^*{}'}$ 
is obtained by slightly deforming $\alpha_0'$ in $\mathbb{R}^2$.
As is suggested in Figure \ref{fig_25}, it is not hard to see that $r_{\mu^*{}',\nu^*{}'}$ is a cubic 
tangency of type I.
\end{proof}
%%%%%%%%%%%%%%%%%%%%%%%%%
%%%%%%%%%%%%%%%%%%%%%%

\begin{remark}
By using an argument similar to that in Lemma \ref{type_1}, one can also suppose that the cubic tangency $r_{\mu^*,\nu^*}$ is of type II.
Two types of cubic homoclinic tangencies corresponding to ours were defined long time ago, for example, see Gonchenko \cite{Go}.
In \cite{Go}, he showed that main bifurcations of periodic points associating with homoclinic cubic 
tangencies of types I and II produce essentially different dynamical phenomena.
\end{remark}

\subsection{Renormalization near generic cubic homoclinic tangencies}\label{SS_renormal}

Let $\{\psi_{\mu,\nu}\}$ be a two-parameter family of $C^r$ diffeomorphisms on $\mathbb{R}^2$ with $\psi_{0,0}=\psi$ and $C^{r-2}$ depending on $\mu,\nu$ for some $r$ with $4\leq r\leq \omega$.
Suppose that  
$\psi=\psi_{0,0}$ has a dissipative saddle fixed point $p$ and a cubic homoclinic tangency $r$ 
of $W^s(p)$ and $W^u(p)$ which unfolds generically with respect to $\{\psi_{\mu,\nu}\}$.
Here the saddle point $p$ is said to be \emph{dissipative} if the eigenvalues $\lambda,\gamma$ of 
the differential $D\psi(p)$ satisfy 
\begin{equation}\label{d_saddle}
0<|\lambda|<1<|\gamma|, \quad |\lambda\gamma|<1.
\end{equation}
The continuation $\{p_{\mu,\nu}\}$ of saddle fixed points of $\psi_{\mu,\nu}$ with $p_{0,0}=p$ is 
well defined for any $(\mu,\nu)$ near $(0,0)$ and the eigenvalues $\lambda=\lambda_{\mu,\nu},\gamma=\gamma_{\mu,\nu}$ of $D\psi_{\mu,\nu}(p_{\mu,\nu})$ also satisfy 
(\ref{d_saddle}).
One can take a coordinate neighborhood $U$ of $p_{\mu,\nu}$ satisfying the following conditions.
\begin{enumerate}[(i)]
\item
The $x$-axis in $U$ is $W_{\mathrm{loc}}^s(p_{\mu,\nu})$ and the $y$-axis is 
$W_{\mathrm{loc}}^u(p_{\mu,\nu})$.
\item
$r=(1,0)$ and $\psi^{-N}(r)=(0,1)$ for some integer $N>0$.
\end{enumerate}

The following lemma is a special case of established results on the renormalization of the first return 
maps near tangencies of arbitrary order, see \cite[Lemma 6]{GST7} or \cite[Lemma 2]{GST3} for 
details.

\begin{lemma}\label{G_Lemma_6}
With the notation as above, for any $k\in \mathbb{N}$, 
there exists a $C^{r-1}$ coordinate transformation 
$(x,y)\mapsto (X,Y)$ and a $C^{r-3}$ parameter transformation $(\mu,\nu)\mapsto (\tilde \mu,\tilde\nu)$ 
such that the first return map $(X,Y)\mapsto (\bar X,\bar Y)$ defined by $\psi_{\mu,\nu}^{N+k}$ near $r$ has the form: 
\begin{align*}
\bar X&= Y+\bar H_{k;1}(\tilde\mu,\tilde\nu,X,Y),\\
\bar Y&= a(\lambda\gamma)^kX+\tilde\nu+\tilde\mu Y+\sigma Y^3+\bar H_{k;2}(\tilde\mu,\tilde\nu,X,Y),
\end{align*}
where $a$ is a non-zero constant, $\sigma$ is $-1$ if the tangency $r$ is of type I and otherwise $1$, and the terms 
$\bar H_{k;j}(\tilde\mu,\tilde\nu,X,Y)$ $(j=1,2)$ are functions 
satisfying the following conditions.
\begin{enumerate}[\rm (i)]
\item
$\bar H_{k;j}$ is well defined on the $R_k$-ball in the $\tilde\mu\tilde\nu XY$-space centered at 
$(0,0,0,0)$ with $\lim_{k\rightarrow \infty}R_k=\infty$.
\item
$\bar H_{k;j}$ is a $C^{r-3}$ function on $(\tilde \mu,\tilde\nu)$ and $C^{r-1}$ on $(X,Y)$.
\item
For any compact subset $K$ of the $\tilde\mu\tilde\nu XY$-space, 
$\Vert \bar H_{k;j}\Vert_{C^3,C^1;K}=O(k\gamma^{-k/2})_{k\rightarrow \infty}$.
\end{enumerate}
\end{lemma}

Here $\Vert \bar H_{k;j}\Vert_{C^3,C^1;K}=O(k\gamma^{-k/2})_{k\rightarrow \infty}$ means that
$\bar H_{k;j}$ itself and all the derivatives of $\bar H_{k;j}$ up to the third order with respect to $X$, $Y$ and 
up to the first order with respect to $\tilde \mu$, $\tilde \nu$ are uniform 
$O(k\gamma^{-k/2})_{k\rightarrow \infty}$-functions on $K$.
This implies that, for any constant $\xi$ with $|\gamma^{-1/2}|<\xi<1$, $\Vert \bar H_{k;j}\Vert_{C^3,C^1;K}=O(\xi^k)_{k\rightarrow \infty}$.

Just replacing \cite[Lemmas 2.1 and 6.6]{KS} by 
the above lemma, we have the following stronger version of \cite[Theorems A and B]{KS}.
In the present theorem, we do not need the local linearizability condition on $\psi$ near saddle fixed points, which was 
crucial in \cite{KS}.

\begin{theorem}\label{T_KS}
Let $\{\psi_{\mu,\nu}\}$ be a two-parameter family in $\mathrm{Diff}^r(\mathbb{R}^2)$ with 
$4\leq r\leq \omega$ 
such that $\psi=\psi_{0,0}$ has a dissipative saddle fixed point $p$ associated with a cubic homoclinic tangency $r$ of type I which unfolds generically with respect to $\{\psi_{\mu,\nu}\}$.
There exist subsets $\mathcal{O}$ and $\mathcal{Z}$ in the $\mu\nu$-space such that the subfamilies $\{\psi_{\mu,\nu}\}_{(\mu,\nu)\in \mathcal{O}}$ and  $\{\psi_{\mu,\nu}\}_{(\mu,\nu)\in \mathcal{Z}}$ satisfy the following conditions.
\begin{enumerate}[\rm (i)]
\item 
$\mathcal{O}$ is an open set with $\mathrm{Cl}(\mathcal{O})\ni (0,0)$.
For any $(\mu,\nu)\in \mathcal{O}$ and a sufficiently small $\varepsilon>0$, there exists a regular curve $c:(-\varepsilon,\varepsilon)\rightarrow \mathcal{O}$ with $c(0)=(\mu,\nu)$ such that the one-parameter family $\{\psi_{c(t)}\}$ exhibits persistent antimonotonic tangencies.
\item 
For any open neighborhood $U$ of $(0,0)$ in the $\mu\nu$-space, $\mathcal{Z}\cap U$ has positive $2$-dimensional Lebesgue measure.
For any $(\mu,\nu)\in \mathcal{Z}$, there exists an integer $n+N>0$ such that $\psi_{\mu,\nu}^{n+N}$ exhibits a cubic polynomial-like strange attractor supported by an SRB measure.
\end{enumerate}
\end{theorem}

\section{Existence of generic cubic tangencies in the H\'enon family}\label{S_3}
In this section, we give the proof of Theorem \ref{main_a}.

\subsection{Saddle fixed points of H\'enon maps}\label{SS_H}

As in Introduction, the H\'enon map $\varphi_{a,b}$ is defined by
$$
\varphi_{a,b}(x,y)=(y, a-bx+y^2).
$$

For any $(a,b)$ close to $(-2,0)$, $\varphi_{a,b}$ has the two fixed points $p_{a,b}^\pm$ with
\begin{equation}\label{p_ab}
p_{a,b}^\pm=(y_{a,b}^\pm,y_{a,b}^\pm),\quad\mbox{where}\quad
y_{a,b}^\pm=\frac{1+b\pm\sqrt{(1+b)^2-4a}}{2}.
\end{equation}
Then the eigenvalues of the differentials $D\varphi_{a,b}(p_{a,b}^\pm)$ are
\begin{equation}\label{e_v}
\lambda_{a,b}^\pm={y_{a,b}^\pm \mp\sqrt{\bigl(y_{a,b}^\pm\bigr)^2-b}},\ \gamma_{a,b}^\pm={y_{a,b}^\pm \pm\sqrt{\bigl(y_{a,b}^\pm\bigr)^2-b}}.
\end{equation}
Since $(\lambda_{a,b}^+,\gamma_{a,b}^+)\rightarrow (0,4)$ and $(\lambda_{a,b}^-,\gamma_{a,b}^-)\rightarrow (0,-2)$ as $(a,b)\rightarrow (-2,0)$, the eigenvalues satisfy
\begin{equation}\label{H_dissipative}
0<|\lambda_{a,b}^\pm| <1<|\gamma_{a,b}^\pm|\quad\mbox{and}\quad | \lambda_{a,b}^\pm \gamma_{a,b}^\pm| <1
\end{equation}
for any $(a,b)$ near $(-2,0)$ with $b\neq 0$.
Thus, both the saddle fixed points $p_{a,b}^\pm$ are dissipative.

\subsection{Outline of proof of Theorem \ref{main_a}}\label{S_outline_b}
Throughout the remainder of this section, $A\sim B$ (resp.\ $\boldsymbol{v}\sim \boldsymbol{w}$) for two real numbers (resp.\ 
vectors) means that one can suppose that $|A-B|<\varepsilon$ (resp. $\Vert\boldsymbol{v}-\boldsymbol{w}\Vert<\varepsilon$) 
for any given $\varepsilon>0$.
We note that $A\sim B$ does not necessarily imply that $A/B$ is close to $1$, e.g.\ $A=\varepsilon/10$ and $B=\varepsilon/1000$.

\begin{figure}[hbt]
\begin{center}
\scalebox{0.8}{\includegraphics[clip]{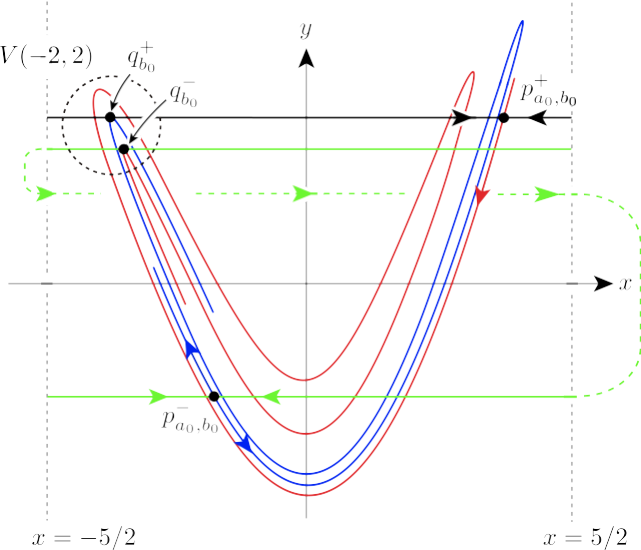}}
\caption{}
\label{fig_31}
\end{center}
\end{figure}

The H\'enon family $\{\varphi_{a,b}\}$ is a $2$-parameter family in $\mathrm{Diff}^\omega(\mathbb{R}^2)$ 
and hence naturally in $\mathrm{Diff}^4(\mathbb{R}^2)$.
Thus, for the proof of Theorem \ref{main_a}, it suffices for us to work in the $C^4$ category. 
By using some results in Kiriki-Li-Soma \cite[Section 3]{KLS}, we have a $C^4$ function $h:I_\varepsilon=(-\varepsilon,\varepsilon)\longrightarrow \mathbb{R}$ with $h(0)=-2$ such that $\varphi_{h(b),b}$ admits a quadratic heteroclinic tangency $q^+_b$ associated with $p^\pm_{h(b),b}$ and contained in a small neighborhood $V(-2,2)$ of $(-2,2)\in \mathbb{R}^2$.
One can also prove that the tangency $q^+_b$ unfolds generically with respect to the $a$-parameter family  $\{\varphi_{a,b(\mathrm{fixed})}\}$.
We will show that there exists $b_0>0$ arbitrarily near $0$ such that $\varphi_{h(b_0),b_0}$ admits a quadratic heteroclinic tangency $q_{b_0}^-$ in $V(-2,2)$ associated with $p_{h(b_0),b_0}^\pm$ such that $q_{b_0}^+,q_{b_0}^-$ are cyclically associated with $p_{h(b_0),b_0}^\pm$ as in Lemma \ref{G_Lemma_5}.
The situation in the present case is illustrated in Figure \ref{fig_31} (cf.\ Figure \ref{fig_21}).
%%%%%%%%%%%%%%%%%%%%%%%%%
%%%%%%%%%%%%%%%%%%%%%%
It remains to show that $q_{b_0}^-$ unfolds generically with respect to the $b$-parameter family $\{\varphi_{h(b),b}\}$.
Unfortunately, the authors do not have any efficient evaluation for the velocity vectors $\boldsymbol{v}^{u,\perp}_b(q_{b_0}^-),\boldsymbol{v}^{s,\perp}_b(q_{b_0}^-)$ defined as in Subsection \ref{S_gqt}.
However, we can get an approximation of the relative velocity vector such that
$$\boldsymbol{v}^{u,\perp}_b(q_{b_0}^-)-\boldsymbol{v}^{s,\perp}_b(q_{b_0}^-)\sim (0,-8\sqrt{2}).$$
From this fact together with (\ref{qt_unfold_3}), we know that $q_{b_0}^-$ unfolds generically.
Then the proof of Theorem \ref{main_a} is completed by taking the new parameter $(\mu,\nu)$ with $\mu=a-h(b), \nu=b-b_0$ and applying Lemma \ref{G_Lemma_5}.

\subsection{Existence of pairs of generically unfolding quadratic tangencies}\label{SS_pair}

In this subsection, we will find a smooth function $a=h(b)$ to find a pair of generically unfolding 
quadratic tangencies associated with $p_{h(b),b}^\pm$ for some $b\in (0,\varepsilon)$.

When $b=0$, $\varphi_{a,0}$ is not a diffeomorphism.
Even in this case, one can define the stable and unstable manifolds associated with $p_{a,0}^+$ in a usual manner.
The stable manifold $W^s(p_{a,0}^+)$ of $\varphi_{a,0}$ contains the horizontal segment $S_{a,0}^+=\{(x,y_{a,0});\vert x\vert\leq 5/2\}$ passing through $p_{a,0}^+$.
By using the Stable Manifold Theorem (for example see Robinson \cite[Chapter 5, Theorem 10.1]{Ro}), 
one can show that, for any $(a,b)$ near $(-2,0)$ (possibly $b=0$),
there exists an almost horizontal segment $S_{a,b}^+\subset W^s(p_{a,b}^+)$ containing $p_{a,b}^+$, 
connecting the vertical lines $x=\pm 5/2$ and $C^4$ depending on $(a,b)$. 
In particular, each $S_{a,b}^+$ has the form
$$S_{a,b}^+=\{(x,\eta_{a,b}^+(x));\vert x\vert \leq 5/2\},$$
where $\eta_{a,b}^+$ is a $C^4$ function $C^4$ depending on $(a,b)$, and the family $\{\eta_{a,b}^+\}$ $C^4$ converges to the constant function $\eta_{a_0,0}^+$ uniformly as $(a,b)\rightarrow (a_0,0)$.

From the definition, the unstable manifold $W^u(p_{a,0}^+)$ consists of the points $q\in \mathbb{R}^2$ which admits a sequence $\{q_n\}_{n=0}^\infty$ in $\mathbb{R}^2$ with $q_0=q$, $q_n\in \varphi_{a,0}^{-1}(q_{n-1})$ for $n=1,2,\dots$ and $\lim_{n\rightarrow \infty}q_n=p_{a,0}^+$.
In particular, $W^u(p_{a,0}^+)$ is contained in the parabolic curve $\mathrm{Im}(\varphi_{a,0})=\{(x,x^2+a); -\infty <x<\infty\}$.
It is not hard to show that
$$W^u(p_{a,0}^+)=\{(x,x^2+a); a \leq x<\infty\}$$
for any $a$ near $-2$.

Let $V(-2,2)$ be a fixed small neighborhood of $(-2,2)$ in the $xy$-plane.
Since $S_{a,0}^+$ is the horizontal line $y_{a,0}=(1+\sqrt{1-4a})/2$, for any $(a,b)$ near $(-2,0)$ and any point $r$ in $S_{a,b}^+$, $\boldsymbol{v}_{a}^{s,\perp}(r)$ is arbitrarily and uniformly close to $\partial y_{a,0}/\partial a|_{a=-2}=-1/3$.
Recall that $\boldsymbol{v}_{a}^{s,\perp}(r)$ is the velocity vector $dr^s(a)/da$ at $r$ defined as in Subsection \ref{S_gqt}, 
where $r^s(a)$ is the intersection point of $S_{a,b(\mathrm{fixed})}^+$ and a short vertical segment passing through $r$.
Thus we have 
\begin{equation}\label{eqn_v_a}
\boldsymbol{v}_{a}^{s,\perp}(r)\sim \Bigl(0,\frac{\partial y_{a,0}}{\partial a}\Big|_{a=-2}\Bigr)=\Bigl(0,-\frac{1}{3}\Bigr).
\end{equation}
Note that probably $\boldsymbol{v}_b^{s,\perp}(r)$ is not constant on $r\in S_{a,0}^+$ though 
$S_{a,0}^+$ is a horizontal segment.
However, if we take $V(-2,2)$ is sufficiently small and $(a,b)$ is sufficiently near $(-2,0)$, then
\begin{equation}\label{eqn_v_a2}
\boldsymbol{v}_{b}^{s,\perp}(\hat r)\sim \frac{d\hat r^s}{db}(0)=:\boldsymbol{v}_0
\end{equation}
for any $\hat r\in S_{a,b}^+\cap V(-2,2)$, where $\hat  r^s(b)$ is the intersection point of $S_{-2,b}^+$ and the vertical line $x=-2$.
The vector $\boldsymbol{v}_0$ will not be evaluated in our proof but it is cancelled out in the 
approximation (\ref{eqn_approx}).

\begin{figure}[hbt]
\begin{center}
\scalebox{0.9}{\includegraphics[clip]{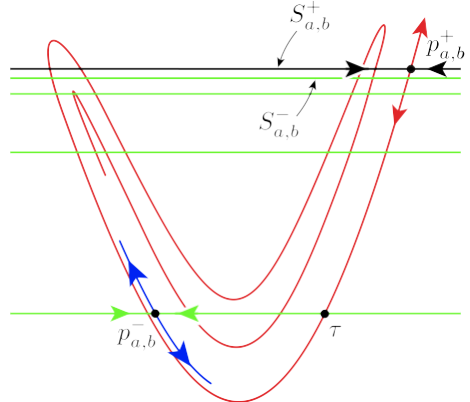}}
\caption{}
\label{fig_32}
\end{center}
\end{figure}

Since $W^s(p_{a,0}^-)$ contains the horizontal line passing through $p_{a,0}^-$, $W^s(p_{a,b}^-)$ and $W^u(p_{a,b}^+)$ have a transverse point $\tau$ for any $(a,b)$ near $(-2,0)$ as illustrated in Figure \ref{fig_32}.
%%%%%%%%%%%%%%%%%%%%%%%%%
%%%%%%%%%%%%%%%%%%%%%%
By the Inclination Lemma, there exists a sequence of curves $S_{a,b;n}^-$ in $W^s(p_{a,b}^-)$ $C^4$ converges to $S_{a,b}^+$ for any $(a,b)$ near $(-2,0)$.
By the Accompanying Lemma \cite[Lemma 4.1]{KS} (see also \cite[Lemma 2.2]{KLS}) together with 
(\ref{eqn_v_a}), (\ref{eqn_v_a2}), one can choose the curve $S_{a,b}^-:=S_{a,b;n}^-$ so that
\begin{equation}\label{v_simeq^-}
\boldsymbol{v}_{a}^{s,\perp}(r')\sim \Bigl(0,-\frac{1}{3}\Bigr),\qquad \boldsymbol{v}_{b}^{s,\perp}({\hat r}')\sim \boldsymbol{v}_0
\end{equation}
for any $r'$ in $S_{a,b}^-$ and ${\hat r}'$ in $S_{a,b}^-\cap V(-2,2)$.

\begin{figure}[htb]
\begin{center}
\scalebox{0.9}{\includegraphics[clip]{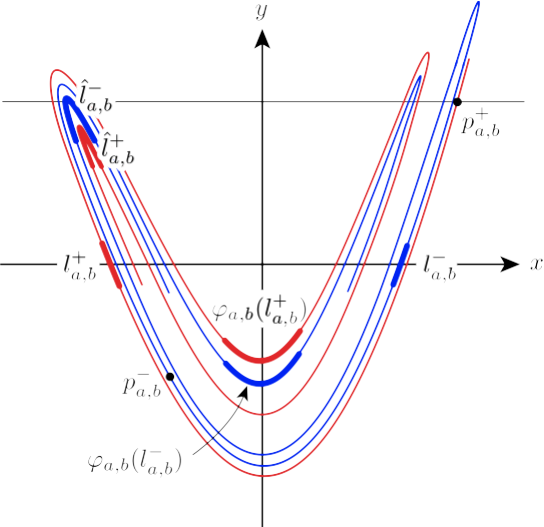}}
\caption{$\varphi_{a,b}$ exchanges the two components of $W_{\mathrm{loc}}^u(p_{a,b}^-)\setminus \{p_{a,b}^-\}$.}
\label{fig_33}
\end{center}
\end{figure}

Let $l_{a,b}^+$ (resp.\ $l_{a,b}^-$) be a short curve in $W^u(p_{a,b}^+)$ (resp.\ $W^u(p_{a,b}^-)$) as illustrated in Figure \ref{fig_33} such that both $\mathrm{Int}(l_{a,b}^\pm)$ meet the $x$-axis transversely.
%%%%%%%%%%%%%%%%%%%%%%%%%
%%%%%%%%%%%%%%%%%%%%%%
Set $\hat l_{a,b}^\pm=\varphi_{a,b}^2(l_{a,b}^\pm)$.
The curve $l_{a,b}^\pm$ is parametrized as $(x_{a,b}^\pm(t),t)$ for any $t$ near $0$, where $x_{a,b}^\pm$ is a $C^4$ function converging uniformly to $x_{-2,0}^\pm(t)=\mp\sqrt{t+2}$ as $(a,b)\rightarrow (-2,0)$.
For $\sigma=\pm$, $\hat l_{a,b}^\sigma$ is parametrized as $\hat r_{a,b}^\sigma (t)=(\hat x_{a,b}^\sigma(t),\hat y_{a,b}^\sigma(t))$, where
\begin{align*}
\hat x_{a,b}^\sigma(t)&=t^2-bx_{a,b}^\sigma(t)+a,\\
\hat y_{a,b}^\sigma(t)&=(\hat x_{a,b}^\sigma(t))^2-bt+a.
\end{align*}
Since $\hat y_{a,b}^\sigma(t)$ $C^4$ converges to $\hat y_{a,0}^\sigma(t)=(t^2+a)^2+a$ as $b\rightarrow 0$, $(\partial/\partial a)\hat y_{a,b}^\sigma(t)$ $C^3$ converges to $(\partial/\partial a)\hat y_{a,0}^\sigma(t)=2(t^2+a)+1$ as $b\rightarrow 0$.
This implies that, for any $(a,b)$ near $(-2,0)$ and any $t$ near $0$, we have
$(\partial/\partial a)\hat y_{a,b}^\sigma(t)\sim -3$ and hence
\begin{equation}\label{v_q}
\boldsymbol{v}^{u,\perp}_a (\hat r_{a,b}^\sigma(t))\sim (0,-3).
\end{equation}

The following lemma asserts that there exists $(a_0,b_0)$ arbitrarily close to $(-2,0)$ such that 
$\varphi_{a_0,b_0}$ has quadratic heteroclinic tangencies $q_{b_0}^\pm$ cyclically associated with $p_{a_0,b_0}^\pm$ which satisfy the assumptions of Lemma \ref{G_Lemma_5} except for 
the generic unfolding property of $q_{b_0}^-$ with respect to the $b$-parameter family $\{\varphi_{a_0,b}\}$.

\begin{lemma}\label{l_bb}
There exists a $C^4$ function $h(b)$ with $h(0)=-2$ defined for any $b$ near $0$ and satisfying the following conditions.
\begin{enumerate}[\rm (i)]
\item
For any non-zero $b$ near $0$, there is a continuation $q_b^+$ of quadratic tangencies of $S^+_{h(b),b}$ and $\hat l_{h(b),b}^-$ each of which unfolds generically with respect to the $a$-parameter family $\{\varphi_{a,b(\mathrm{fixed})}\}$.
\item
For any sufficiently small $b_1>0$, there is $b_0$ with $0<b_0<b_1$ and an arc $S_{h(b_0),b_0}^-$ in $W^s(p_{h(b_0),b_0 }^-)$ such that $S^-_{h(b_0),b_0}$ and $\hat l_{h(b_0),b_0}^+$ have a quadratic tangency $q_{b_0}^-$.
\end{enumerate}
\end{lemma}
\begin{proof}
By an argument quite similar to that in \cite[Section 3]{KLS}, there exists a $C^4$ function $h:I_\varepsilon=(-\varepsilon,\varepsilon)\longrightarrow \mathbb{R}$ for a sufficiently small $\varepsilon>0$ with $h(0)=-2$ and such that $S^+_{ h(b),b}$ and $\hat l_{h(b),b}^-$ have a quadratic tangency $q_{b}^+$ which is a continuation on $b\in (-\varepsilon,\varepsilon)$.
Note that, in \cite{KLS}, the tangency is a homoclinic one associated with $p_{a,b}^+$, but in the present case, 
the tangency is a heteroclinic one associated with $p_{a,b}^\pm$.
However, a similar argument works since the unstable manifold $W^u(p_{a,0}^-)$ has the form:
$$W^u(p_{a,0}^-)=\{(x,x^2+a);\, a\leq x<\infty\}$$
as well as $W^u(p_{a,0}^+)$.

Fix $b_1\in (0,\varepsilon)$ arbitrarily.
If necessary retaking $S_{h(b_1),b_1}^-$ again, we may assume that the level of 
$S_{h(b_1),b_1}^-$ is higher than that of $\hat l_{h(b_1),b_1}^+$, see Figure \ref{fig_34}\,(a).
%%%%%%%%%%%%%%%%%%%%%%%%%
\begin{figure}[hbt]
\begin{center}
\scalebox{0.7}{\includegraphics[clip]{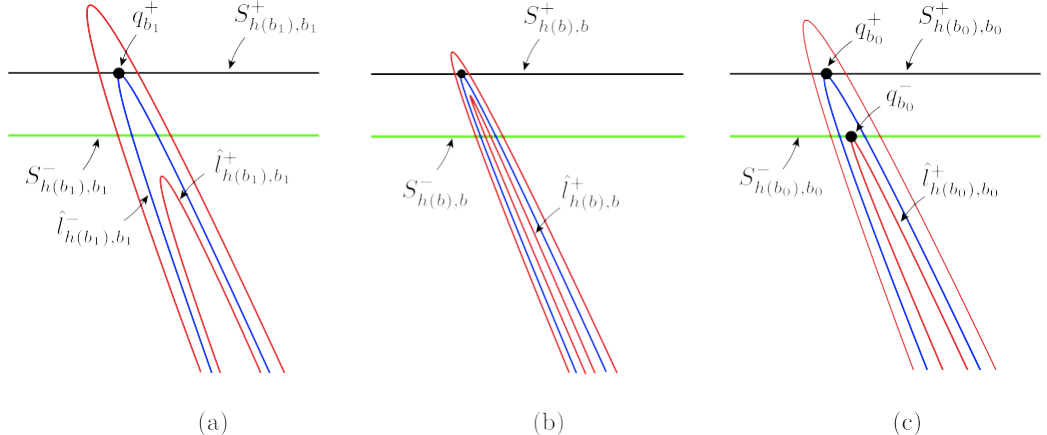}}
\caption{}
\label{fig_34}
\end{center}
\end{figure}
%%%%%%%%%%%%%%%%%%%%%%
Note that the maximal point of $\hat l_{h(b),b}^+$ approaches to $S_{h(b),b}^+$ as $b\searrow 0$, see Figure \ref{fig_34}\,(b).
On the other hand, $S_{h(b),b}^-$ converges to $S_{-2,0}^-$, which is a horizontal segment the level of which is strictly 
lower than that of $S_{-2,0}^+$.
By the Intermediate Value Theorem, there exists $b_0$ with $0<b_0<b_1$ such that $S_{h(b_0),b_0}^-$ and $\hat l_{h(b_0),b_0}$ have a quadratic tangency $q^-_{b_0}$, see Figure \ref{fig_34}\,(c).
The tangency $q_{b_0}^+$ is equal to $\hat r_{h(b_0),b_0}^-(t_{b_0})$ for some $t_{b_0}$ near $0$.
Thus the generic condition (\ref{qt_unfold_3}) and the approximations (\ref{eqn_v_a}), (\ref{v_q}) show that $q_{b_0}^+$ unfolds generically with respect to $\{\varphi_{a,b(\mathrm{fixed})}\}$.
This completes the proof.
\end{proof}

For the completion of the proof of Theorem \ref{main_a}, it suffices to show that the tangency $q_{b_0}^-=\hat r_{h(b_0),b_0}^-$ given in Lemma \ref{l_bb} unfolds generically with respect to the $b$-parameter family $\{\varphi_{h(b),b}\}$.
Then $\{\varphi_{a,b}\}$ satisfies the conditions (i) and (ii) of Lemma 
\ref{G_Lemma_5} with respect to the new parameters
$$\mu=a-h(b),\quad \nu=b-b_0.$$
We denote the graph function of $S_{a,b}^-$ given in Lemma \ref{l_bb} by $\eta_{a,b}^-$ and the subscription pair `$h(b),b$' only by `$b$', e.g.\ $\eta_{h(b),b}^\pm=\eta_b^\pm$.
Recall that $\eta^+_b$ is the graph function of $S_{h(b),b}^+$ and note that $\eta_0^\pm(x)$ $(|x|\leq 5/2)$ are constants, denote them by $c^\pm$.

\begin{proof}[Proof of Theorem \ref{main_a}]

For $\sigma=\pm$, we set $-\sigma=\mp$.
Recall that $\hat l_b^\sigma$ has the regular curve  parametrization $\bigl(\hat x_b^\sigma(t),(\hat x_b^\sigma(t))^2-bt+h(b)\bigr)$,
where 
$$
\hat x_b^\sigma(t)=t^2-bx_b^\sigma(t)+h(b).
$$
Then the function $\theta_b^\sigma(t)$ of $t$ corresponding to (\ref{theta_mu}) is defined by
\begin{equation}\label{theta_b}
\theta_b^\sigma(t)=(\hat x_b^\sigma(t))^2-bt+h(b)-\eta_b^{-\sigma}(\hat x_b^\sigma(t)).
\end{equation}

Since $\theta_b^\sigma(t)$ $C^4$ converges to $\theta_0^\sigma(t)=(t^2-2)^2+c^\sigma$ as $b\rightarrow 
0$, 
it follows that $\dot\theta_b^\sigma(t)$  $C^3$ converges  to $\dot\theta_0^\sigma(t)=4t^3-8t$ and $\ddot\theta_b^\sigma(t)$ $C^2$ converges to 
$\ddot\theta_0^\sigma(t)=12t^2-8\sim -8\neq 0$ as $b\rightarrow 
0$, where the `dot' represents the derivative of a function by $t$.
From this fact, we have a $C^3$ function $t^\sigma_b$ of $b$ with $t^\sigma_0=0$ and $\dot\theta_b^\sigma(t_b^\sigma)=0$.
Note that, since $q_b^+$ is in $S_b^+$ for any $b$ near $0$, $\theta_b^-(t_b^-)$ is a zero-constant function of $b$.

We set $x_b^\sigma(t)=x^\sigma(b,t)$ and $\eta_b^\sigma(x)=\eta^\sigma(b,x)$, and suppose that the `prime' represents the derivative of a function by $b$.
For example, ${x_b^\sigma}'(t)=(\partial/\partial b) x^\sigma(b,t)$ and ${\eta_b^\sigma}'(x)=(\partial/\partial b) \eta^\sigma(b,x)$.

Now, we show that, for any sufficiently small $b>0$,
\begin{equation}\label{d_theta_b}
\frac{d}{db}\theta_b^+(t_b^+)\sim -8\sqrt{2}.
\end{equation}
Set $t_b^+=t_b$, $x_b^+=x_b$ and $\hat x_b^+=\hat x_b$ for short.
By (\ref{theta_b}),
\begin{align*}
\frac{d\theta_b^+(t_b)}{db}&=2\hat x_b(t_b)\bigl(2t_bt_b'-x_b(t_b)-bx_b'(t_b)-b\dot x_b(t_b)t_b'+h'(b)\bigr)\\
&\qquad\qquad -t_b-bt_b'+h'(b)-\eta_b^-{}'(\hat x_b(t_b))-\frac{\partial \eta^-}{\partial x}(b,\hat x_b(t_b))(\hat x_b(t_b))'.
\end{align*}
Since $\eta_0^-$ is a constant function of $x$, we have $(\partial/\partial x)\eta_0^-=0$.
Since moreover $t_0=0$, $x_0(0)=-\sqrt{2}$ and  $\hat x_0(0)=h(0)\sim -2$, it follows that
$$
\frac{d\theta_b^+(t_b^+)}{db}\sim -4\sqrt{2}-3h'(0)-\eta^-{}'(0,-2).
$$
Since $\theta_b^-(t_b^-)$ is a zero-constant function and since $x_0^-(0)=\sqrt2$, $\hat x_0^-(0)=h(0)\sim -2$, we have similarly
$$
0=\frac{d\theta_b^-(t_b^-)}{db}\sim 4\sqrt{2}-3h'(0)-\eta^+{}'(0,-2).
$$
By (\ref{eqn_v_a2}) and (\ref{v_simeq^-}),
$$\bigl(0,\eta^-{}'(0,-2)\bigr)\sim \boldsymbol{v}_0,\quad \bigl(0,\eta^+{}'(0,-2)\bigr)\sim \boldsymbol{v}_0.$$
The above four approximations imply (\ref{d_theta_b}).

From (\ref{qt_unfold_3}) and (\ref{d_theta_b}), one can get our desired approximation
\begin{equation}\label{eqn_approx}
\boldsymbol{v}^{u,\perp}_b(q_{b_0}^-)-\boldsymbol{v}^{s,\perp}_b(q_{b_0}^-)= \Bigl(0,\frac{d\theta_b^+(t_b^+)}{db}(b_0)\Bigr)\sim (0,-8\sqrt{2})
\end{equation}
by taking the $b_0>0$ in Lemma \ref{l_bb} sufficiently small.
It follows that $q_{b_0}^-$ unfolds generically with respect to the $b$-parameter family $\{\varphi_{h(b),b}\}$.

By Lemma \ref{G_Lemma_5}, there exists a cubic tangency $r_{a^*,b^*}$ associated with $p_{a^*,b^*}^+$ for 
some $(a^*,b^*)$ arbitrarily close to $(h(b_0),b_0)\sim (-2,0)$ and satisfying the properties in Theorem \ref{main_a}.
Since $b_0>0$, by (\ref{e_v}) the both eigenvalues of $D\psi_{h(b_0),b_0}(p_{h(b_0),b_0}^+)$ are positive and 
those of $D\psi_{h(b_0),b_0}(p_{h(b_0),b_0}^-)$ are negative.
Hence, by Lemma \ref{type_1}, one can suppose that the tangency $r_{a^*,b^*}$ is of type I.
\end{proof}

\subsection*{Acknowledgments}
The authors would like to thank the referee for valuable comments and suggestions, 
according to which many parts of this paper are improved and corrected.

\end{document}